\pdfoutput=1
\documentclass{ifacconf}
\usepackage[algo2e, ruled, vlined]{algorithm2e}

\usepackage{graphicx}      
\usepackage{natbib}        

\usepackage{xargs}
\usepackage{amsmath}
\usepackage{amssymb}
\usepackage{mathtools}
\usepackage{xcolor}
\usepackage{float}
\usepackage{mathrsfs}
\usepackage{relsize}
\usepackage{url}

\SetKwFor{While}{while}{}{}%
\SetKwFor{For}{for}{}{}%
\SetKwIF{If}{ElseIf}{Else}{if}{}{else if}{else}{end if}%

\makeatletter
 \let\old@ssect\@ssect  
 \makeatother
\usepackage[colorlinks=true,linkcolor=black,anchorcolor=black,citecolor=black,filecolor=black,menucolor=black,runcolor=black,urlcolor=black]{hyperref}
\makeatletter
 \def\@ssect#1#2#3#4#5#6{
   \NR@gettitle{#6}  
   \old@ssect{#1}{#2}{#3}{#4}{#5}{#6}  
 }
 \makeatother
\renewcommand\#{\protect\scalebox{0.6}{\protect\raisebox{0.4ex}{\char"0023}}}

\hypersetup{pdfauthor={Pieter Pas, Andreas Themelis, and Panagiotis Patrinos},pdftitle={Gauss-Newton meets PANOC: A fast and globally convergent algorithm for nonlinear optimal control}}

\newcommand\Reals{\mathrm{I\!R}}
\newcommand\Nats{\mathrm{I\!N}}

\newcommand\R{\Reals}
\newcommand\N{\Nats}

\newcommand\minimize{\operatorname*{\mathbf{minimize}}}

\newcommand\subjto{\operatorname*{\mathbf{subject\;to}}}
\newcommand\prox{\operatorname{\mathbf{prox}}}

\newcommand\proj[1]{\operatorname{\mathbf{\Pi}_\mathnormal{#1}}}
\newcommand\dist{\operatorname{\mathbf{dist}}}
\newcommand\Id{\mathrm{Id}}
\newcommand\I{\mathrm{I}}

\newcommand\defeq{\triangleq}

\DeclareMathOperator*{\argmin}{arg\,min}

\newcommand\norm[1]{\left\| {#1} \vphantom{X} \right\|}
\newcommand\normsq[1]{\norm{#1}^2}
\newcommand\defset[2]{\left\{ {#1} \;\middle|\; {#2} \right\}}
\newcommand\tp[1]{#1^\top}
\newcommand\ttp[1]{#1^{\!\top\!}}

\newcommand\inv[1]{#1^{-1}}

\newcommand\iddots{\mathinner{
  \kern1mu\raise1pt{.}
  \kern2mu\raise4pt{.}
  \kern2mu\raise7pt{\Rule{0pt}{7pt}{0pt}.}
  \kern1mu
}}

\newcommand\jac{\mathrm{J}}

\newcommand\possdef{\succeq}

\newcommand\delx{\Delta x}
\newcommand\delu{\Delta u}

\newcommand\U{U}
\newcommand\D{D}
\newcommand\K{\mathcal K}
\newcommand\J{\mathcal J}
\newcommand\nnu{{n_u}}
\newcommand\nnx{{n_x}}

\newcommand\lt{<}

\newcommand\fbop{T_\gamma}
\newcommand\fpr{R_\gamma}
\newcommandx\fbe[1][1=\gamma]{\varphi^\mathrm{FB}_{#1}}

\DeclareMathOperator*{\interior}{\mathbf{int}}
\DeclareMathOperator*{\boundary}{\mathbf{bdry}}

\newcommand\Hgn{\hat\nabla^2_{\!\text{GN}}}
\newcommand\egn{\delta^2_{\text{GN}}}
\newcommand\jacF{\jac_{\!F}}
\newcommand\jacf{\jac_{\!f}}
\newcommand\hhbar{\hbar} 
\newcommand\barxuk{\bar x^k\!, \bar u^k}

\newcommand\panoc{PANOC}
\newcommand\Panoc{PANOC}

\newcommand\Cpp{C\nolinebreak[4]\hspace{-.02em}\raisebox{.36ex}{\relsize{-3}{\textbf{++}}}}

\def\mybigtimes{\mathop{\mathchoice{
   \vcenter{\hbox to10bp{\vrule height15bp width0pt \pdfliteral{
   q 1 J .8 w 0 1 m 10 14 l S 0 14 m 10 1 l S Q
}\hss}}}{
   \vcenter{\hbox to10bp{\kern1bp\vrule height10bp width0pt \pdfliteral{
   q 1 J .65 w 0 0 m 8 10 l S 0 10 m 8 0 l S Q
}\hss}}}{\times}{\times}
}}

\begin{document}
\begin{frontmatter}

\title{Gauss--Newton meets PANOC: \\A fast and globally convergent algorithm for nonlinear optimal control\thanksref{footnoteinfo}} 

\thanks[footnoteinfo]{This work is supported by: the Research Foundation Flanders (FWO) PhD grant No. 11M9523N and research projects G081222N, G033822N, G0A0920N; European Union’s Horizon 2020 research and innovation programme under the Marie Skłodowska-Curie grant agreement No. 953348;
Research Council KU Leuven C1 project No. C14/18/068; Fonds de la Recherche Scientifique – FNRS; Fonds Wetenschappelijk Onderzoek – Vlaanderen under EOS project No. 30468160 (SeLMA); and the Japan Society for the Promotion of Science
(JSPS) KAKENHI grant JP21K17710.}

\author[kul]{Pieter Pas} 
\author[kyushu]{Andreas Themelis} 
\author[kul]{Panagiotis Patrinos}

\address[kul]{Department of Electrical Engineering (ESAT--STADIUS), \\
KU Leuven, Kasteelpark Arenberg 10, Leuven, 3001, Belgium \\ (e-mail: \{pieter.pas,panos.patrinos\}@esat.kuleuven.be)}
\address[kyushu]{Faculty of Information Science and Electrical Engineering (ISEE), \\
   Kyushu University, 744 Motooka, Nishi-ku, 819-0395 Fukuoka, Japan \\ (e-mail: andreas.themelis@ees.kyushu-u.ac.jp)}

\begin{abstract}                
    \Panoc{} is an algorithm for nonconvex optimization that has recently gained
    popularity in real-time control applications due to its fast, global convergence.
    The present work proposes a variant of \panoc{} that makes use of
    Gauss--Newton directions to accelerate the method. Furthermore, we show that
    when applied to optimal control problems, the computation of this 
    Gauss--Newton step can be cast as a linear quadratic regulator (LQR) problem,
    allowing for an efficient solution through the Riccati recursion.
    Finally, we demonstrate that the proposed algorithm is more than twice as
    fast as the traditional L--BFGS variant of \panoc{} when applied to an
    optimal control benchmark problem, and that the performance scales favorably
    with increasing horizon length.
\end{abstract}

\end{frontmatter}

\section{Introduction}
The ever increasing scale and complexity of models used in optimal control
applications necessitate the development of efficient numerical solvers for
large-scale, nonconvex optimization. One such solver is \panoc{},
the \textit{Proximal Averaged Newton-type method for Optimality Conditions}
\cite[]{panoc}, which has proven successful in
real-time model predictive control (MPC) applications
\cite[]{8550253,small_AerialNavigationObstructed_2019,9625659}. 
Various implementations are available, in \Cpp{} \cite[]{pas2022alpaqa}, 
Rust \cite[]{open}, and Julia \cite[]{Stella_ProximalAlgorithms_jl_Proximal_algorithms}.
The appeal of an algorithm like \panoc{} is that it enjoys fast convergence
thanks to its Newton-type directions, without giving up any theoretic guarantees
about global convergence \cite[]{de_marchi_proximal_2022}.

In the original \panoc{} publication, the limited-memory BFGS (L--BFGS) method
was used to generate fast Newton-type directions. In \cite[]{pas2022alpaqa},
the structure of box-constrained problems was exploited to apply L--BFGS more
effectively by reducing the Newton system to a lower-dimensional one after
eliminating active constraints. The present work continues the search for
faster and more effective Newton-type directions by exploiting the specific
structure of optimal control problems (OCPs).

The remainder of this paper is structured as follows. 
In Section~\ref{sec:gn-panoc}, we explore a linear Newton approximation (LNA)
of the fixed-point residual mapping that lies at the core of \panoc{}.
By using a Gauss--Newton approximation, the high computational cost of
evaluating second-order derivatives is avoided.
In Section~\ref{sec:optimal-control}, we go on to apply this Gauss--Newton
variant of \panoc{} to an input-constrained, nonconvex optimal control problem,
and show that the computation of the Gauss--Newton step corresponds to the solution of an
equality-constrained linear quadratic regulator (LQR) problem.
Section~\ref{sec:algorithm} covers efficient algorithms for solving this LQR
problem by using the Riccati recursion. Pseudocode for the full algorithm is
provided, as well as a brief discussion of the computational cost of the
operations involved.
The performance of the resulting algorithm is validated in Section~\ref{sec:experiments},
where it is applied to a challenging model predictive control benchmark.
We report a speedup by a factor of two compared to the L--BFGS version of
\panoc{}. Finally, Section~\ref{sec:conclusion} concludes with a recapitulation
of the main results and a discussion of future work.

\subsection{Notation}
Let \([a,b]\) denote the closed interval from \(a\) to \(b\). \(\N_{[i,j]} \defeq [i,j] \cap \N\) is the inclusive range of natural numbers from \(i\) to \(j\).
\(\overline \R \defeq \R \cup \{+\infty\}\) is the set of extended real values.
$x_i$ refers to the $i$'th component of $x \in \R^n$. 
Given an index set $\mathcal{I} = \{i_1,\, \dots,\, i_m\} \subseteq \N_{[1,n]}$,
we use the shorthand $x_{\mathcal{I}} = (x_{i_1},\, \dots,\, x_{i_m})$.
Given a matrix \(A\in \R^{n\times m}\), \(A{\scriptstyle[\mathcal I\!, \mathcal J]} \in \R^{\#\mathcal I \times \#\!\mathcal J}\)
denotes the matrix that consists of all elements of \(A\) with row indices
in index set \(\mathcal I\) and column indices in \(\mathcal J\); a dot is
used to denote all indices, e.g. \(A{\scriptstyle[\mathcal I\!,\,\cdot\,]}\)
selects the complete rows of \(A\) with row indices in \(\mathcal I\).
For \(u,v \in \R^n\), let \(u \le v\) denote the component-wise comparison.
In the context of receding horizon problems, the vector \(u\in \R^{N\nnu}\)
without superscript refers to the concatenation of all vectors
\(u^k \in \R^\nnu\) for each time step \(k\) in the horizon.
Given a positive definite matrix \(R\), define the \(R\)-norm as
\(\|x\|_R \defeq \sqrt{x^{\!\top} R\, x}\); in the absence of a subscript, \(\|x\|\) refers to the Euclidean norm.
The indicator function \(\delta_\U\) of a set \(\U\) is zero if its argument is
an element of \(\U\) and \(+\infty\) otherwise.
The proximal operator of a function \(g : \R^n \rightarrow \overline \R\) is
defined as \(\prox_{g}(x) \defeq 
\argmin_w \big\{ g(w) + \tfrac{1}{2} \left\| w - x \right\|^2 \big\}\),
with as a special case \(\prox_{\delta_{\U}}(x) =
\argmin_{w\in \U} \big\{ \tfrac{1}{2}\left\| w - x \right\|^2 \big\} \defeq \proj{\U}(x)\)

\cite[\S 1.G]{RockafellarVariationalAnalysis}. Denote the distance between a point \(x\) and a closed set \(D\) by \(\dist_D(x) = \|x - \proj{D}(x)\|\). \\
Let \(f : \R^n \to \R^p\) and \(g : \R^m \to \R^q\),
then \((f \times g) : \R^n \times \R^m \to \R^p \times \R^q : (x, y) \mapsto (f(x), g(y))\) is their Cartesian product,
and if \(p = q\), their reduced sum is defined as \((f\oplus g) : \R^n \times \R^m \to \R^p : (x, y) \mapsto f(x) + g(y)\). \\
For a function \(F : \R^n \to \R^m\), denote its Jacobian matrix by \(\jacF : \R^n \to \R^{m \times n}\).
For multivariate functions, a superscript is used to refer to the variable with respect to which to differentiate, 
\(\jacF^x \defeq \frac{\partial F}{\partial x} = \tp\nabla_{\!x}\!F\). The Clarke generalized Jacobian of \(F\)
is denoted by \(\partial_CF\) \cite[]{clarke1990optimization}, and for a differentiable function \(f : \R^n \to \R\),
define the generalized Hessian matrix \(\partial^2f \defeq \partial_C(\nabla f)\).

\section{Gauss--Newton acceleration of PANOC}\label{sec:gn-panoc}
We consider optimization problems of the general form
\begin{equation}\label{eq:P} \tag{P}
    \begin{aligned}
        &&\minimize_u &\quad \psi(u) + g(u), \\
    \end{aligned}
\end{equation}
where \(\psi : \R^n \to \R\) has a locally Lipschitz-continuous gradient but is
not necessarily convex, and where \(g(u) : \R^n \to \overline \R\) is proper,
lower semicontinuous, and
\(\gamma_g\)-prox-bounded, but possibly nonsmooth and nonconvex.
Problems of this form can be tackled
using the proximal gradient method, or accelerated variants thereof, such as the
\panoc{} algorithm \cite[]{panoc,de_marchi_proximal_2022}.

\subsection{Linear Newton approximations for \panoc}

Local solutions to \eqref{eq:P} correspond to fixed points of the \textit{forward-backward operator}
\(\fbop(u) \defeq \prox_{\gamma g}\big(u - \gamma\nabla\psi(u)\big)\),
and are characterized by the nonlinear inclusion 
\(0 \in \fpr(u)\), 
where \(\fpr \defeq \inv\gamma(\Id - \fbop)\) is the \textit{fixed-point residual} of \(\fbop\).
Traditionally, \panoc{} applies the L--BFGS quasi-Newton method to this root-finding
problem to achieve fast convergence. A line search over the \textit{forward-backward envelope} \(\fbe\)
is used as a globalization strategy.

This paper explores alternative directions to accelerate \panoc{} by studying
generalized Jacobians to construct a \textit{linear Newton approximation} (LNA) \cite[]{pang}
of the fixed-point residual \(\fpr\).

\newcommand\proxjac{B}
\begin{prop}\label{prop:lna} (LNA scheme for \(\fpr\))\\
Suppose that \(\nabla\psi\) is semismooth around \(\bar u\in \R^n\) and
that \(\prox_{\gamma g}\) with \(\gamma > 0\) is semismooth at
\(\bar u - \gamma\nabla \psi(\bar u)\). Then,
\begin{equation}\label{eq:lna}
    H_\gamma(u) \defeq \inv\gamma \I - \proxjac(u)\big(\inv\gamma\I - \partial^2\psi(u)\big),
\end{equation}
where \(\proxjac(u) = \partial_C \prox_{\gamma g} \big(u - \gamma\nabla\psi(u)\big)\)
and \(\partial^2\psi(u) = \partial_C\big(\nabla\psi(u)\big)\),
furnishes an LNA scheme for \(\fpr\) at \(\bar u\).
\cite[Lem.~6]{6760233}
\cite[Prop.~3.7]{fbtruncnewton}
\cite[\S 15.4.13]{themelis2019acceleration}
\end{prop}
\begin{pf}
    Because of the semismoothness of \(\prox_{\gamma g}\) and \(\nabla\psi\),
    \(\proxjac(u)\) is an LNA scheme for \(\prox_{\gamma g}\) at \(\bar u - \gamma\nabla\psi(\bar u)\),
    and \(\I - \gamma\partial^2\psi(u) = \partial_C\big( u - \gamma\nabla\psi(u)\big)\) is 
    an LNA scheme for \(\Id - \gamma\nabla\psi\) at \(\bar u\). By \cite[Thm.~7.5.17]{pang},
    the product \(\proxjac(u)\big(\inv\gamma\I - \partial^2\psi(u)\big)\) is an LNA scheme for
    the composition \(\fbop = \prox_{\gamma g} \circ (\Id - \gamma\nabla\psi)\) at \(\bar u\).
    \hfill\(\square\)
\end{pf}

This proposition motivates using a solution \(\delu\) of the Newton system \(H_\gamma(\bar u)\,\delu = -\fpr(\bar u)\)
as an update direction for \panoc{}, using the LNA around the current iterate
\(\bar u\).

\subsection{Structured PANOC}
In the case where the nonsmooth term \(g\) in \eqref{eq:P} is the indicator of a closed rectangular
box \(\U\), i.e. \(g \defeq \delta_\U\),
\(\prox_g\) is a separable projection. This structure can be exploited to reduce
the dimension of the Newton system \cite[\S III]{pas2022alpaqa}.

Represent the box \(U \defeq \mybigtimes_{i=1}^n U_i\)
as a Cartesian product of one-dimensional intervals. Then, \(\proxjac(u) = \partial_C \proj{\U}\!\big(u - \gamma\nabla\psi(u)\big)\) is a
set of diagonal matrices with
\begin{equation}
    \proxjac(u)_{ii} \in \begin{cases}
        \{0\} &\text{ if } u_i - \gamma\nabla_{\!i}\psi(u) \not\in \U_i, \\
        \{1\} &\text{ if } u_i - \gamma\nabla_{\!i}\psi(u) \in \interior \U_i, \\
        [0, 1] &\text{ if } u_i - \gamma\nabla_{\!i}\psi(u) \in \boundary \U_i. \\
    \end{cases}
\end{equation}
Motivated by these different cases, let us define the index sets \(\K(u) \defeq \defset{i\in \N_{[1,\,n]}}{u_i - \gamma \nabla_{\!i} \psi(u) \not\in \interior \U_i}\)
and \(\J(u) \defeq \defset{i\in \N_{[1,\,n]}}{u_i - \gamma \nabla_{\!i} \psi(u) \in \interior \U_i}\)
of active and inactive constraints respectively,
and choose \(\hat \proxjac(u) \in \proxjac(u)\), defining \(\hat \proxjac(u)_{ii} \defeq 0\) if \(i \in \K(u)\) and \(\hat \proxjac(u)_{ii} \defeq 1\) if \(i \in \J(u)\).

\newcommand\newtdir{\delu}
By permutation of \eqref{eq:lna}, the Newton step \(\newtdir\) at a point \(\bar u\) can then be computed by solving
the system
\begin{equation}\label{eq:struc-panoc-sys}
    \left\{
    \begin{aligned}
        \newtdir_\K &= \bar u_\K - \fbop(\bar u)_\K, \\
        \partial^2_{\!\J\!\J}\psi(\bar u)\, \newtdir_{\!\J} &= -\nabla_{\!\!\J}\psi(\bar u) - \partial^2_{\!\J\!\K} \psi(\bar u)\, \newtdir_\K. \\
    \end{aligned}\right.
\end{equation}

\begin{figure*}[t]
\normalsize
\setcounter{equation}{5}
\begin{align}
    \label{eq:ocp-qp}\tag{P-ELQR}
    &\begin{aligned}
        &\minimize_{\delx,\delu}&&\quad \tfrac12 \sum_{k=0}^{N-1} 
        \ttp{\begin{pmatrix} \delx^k \\ \delu^k \end{pmatrix}}\! \begin{pmatrix}
            Q_k & \tp S_k \\ S_k & R_k
        \end{pmatrix} \begin{pmatrix} \delx^k \\ \delu^k \end{pmatrix} 
        + \tfrac12 \left( \delx^{N} \right)^\top\! Q_N \left( \delx^{N} \right)
        + \sum_{k=0}^{N-1} \ttp{\begin{pmatrix} q^k \\ r^k \end{pmatrix}}\!
        \begin{pmatrix} \delx^k \\ \delu^k \end{pmatrix} 
        + \left(q^{N}\right)^\top\! \left(\delx^N\right) \\
        &\subjto &&\quad \delx^0 = 0 \\
        &&&\quad \delx^{k+1} = \makebox[0pt][l]{$\displaystyle{A_k \delx^k + B_k \delu^k}$}\hphantom{A_k \delx^k + \hat B_k \delu^k_{\!\J} + \hat c_k} \quad\quad \raisebox{0.3ex}{$\scriptstyle(0 \le k \lt N)$} \\
        &&&\quad \delu_\K = u_\K - \fbop(u)_\K \\[1em]
    \end{aligned} \hspace{-2.5em} \\
    \label{eq:ocp-qp-elim}\tag{P-LQR}
    &\begin{aligned}
        &\minimize_{\delx,\delu_{\!\J}}&&\quad \tfrac12 \sum_{k=0}^{N-1} 
        \ttp{\begin{pmatrix} \delx^k \\ \delu^k_{\!\J} \end{pmatrix}}\! \begin{pmatrix}
            Q_k & \tp{\hat{S}}_k \\ \hat{S}_k & \hat{R}_k
        \end{pmatrix} \begin{pmatrix} \delx^k \\ \delu^k_{\!\J} \end{pmatrix} 
        + \tfrac12 \left( \delx^{N} \right)^\top\! Q_N \left( \delx^{N} \right)
        + \sum_{k=0}^{N-1} \ttp{\begin{pmatrix} \hat q^k \\ \hat r^k \end{pmatrix}}\!
        \begin{pmatrix} \delx^k \\ \delu^k_{\!\J} \end{pmatrix} 
        + \left( \hat q^{N} \right)^\top\! \left( \delx^N \right) \\
        &\subjto &&\quad \delx^0 = 0 \\
        &&&\quad \delx^{k+1} = A_k \delx^k + \hat B_k \delu^k_{\!\J} + \hat c_k \quad\quad \raisebox{0.3ex}{$\scriptstyle(0 \le k \lt N)$} \\
    \end{aligned} \hspace{-2.5em}
\end{align}
\hrulefill
\end{figure*}

\subsection{Gauss--Newton approximation}
We will now specialize to problems where the smooth term is a composition
\(\psi(u) \defeq \ell\big(F(u)\big)\)
of \(\ell : \R^m \to \R\) convex and \(F : \R^n \to \R^m\).
Considering the computational cost of evaluating and factorizing the second-order
derivatives of \(\psi\), the proposed method approximates \eqref{eq:struc-panoc-sys}
using the Gauss--Newton matrix \(\Hgn \defeq \tp{\jacF(u)}\, \partial^2 \ell\big(F(u)\big)\, \jacF(u)\)
\cite[\S3]{schraudolph2002fast}.
\begin{rem}
For \(\psi \in C^2\), we have \(\nabla^2\psi = \Hgn + \egn\) with
\(\egn(u) \defeq \sum_{i=1}^m \nabla_{\!i\,} \ell\big(F(u)\big)\, \nabla^2 F_i(u)\).
If the function \(F\) is linear around a solution \(u^\star\), or if \(F(u^\star)\)
is a stationary point of \(\ell\), the error term \(\egn\) vanishes, and the 
Gauss--Newton approximation approaches the true Hessian matrix of \(\psi\).
\end{rem}
Substituting \(\partial^2\psi\) by \(\Hgn\) in \eqref{eq:struc-panoc-sys} and 
writing the solution to the resulting system as the solution of an equality
constrained quadratic program yields
\begin{equation}\label{eq:eqp}\tag{GN-QP}
    \begin{aligned}
            &\minimize_{\delu}&& \tfrac12\, \tp \delu \Hgn(\bar u)\, \delu + \ttp{\nabla\psi(\bar u)} \delu \\
            &\subjto && \delu_\K = u_\K - \fbop(\bar u)_\K.
    \end{aligned}
\end{equation}
The following sections explore methods for efficiently solving this 
Gauss--Newton QP by making use of the particular structure of finite-horizon
optimal control problems. The Gauss--Newton step \(\delu\) can then be used as
an accelerated direction for \panoc.

\section{Optimal control}\label{sec:optimal-control}
This section explores how optimal control problems arising in model predictive
control applications fit into the optimization framework from the previous
section, and how their specific structure can be exploited to compute
Gauss--Newton directions efficiently.

\subsection{Problem formulation}

Consider the following general formulation of a nonlinear optimal control
problem with finite horizon \(N\).
\newcommand\xinit{x_\text{init}}
\newcommand\xref{x_\text{r}}
\newcommand\uref{u_\text{r}}
\begin{equation}\label{eq:OCP} \tag{OCP}\hspace{-0.8em}
    \begin{aligned}
        &\minimize_{u,x} && \sum_{k=0}^{N-1} \ell_k\big(h_k(x^k, u^k)\big) + \ell_N\big(h_N(x^N)\big)\hspace{-0.8em} \\
        &\subjto && u \in \U \\
        &&& x^0 = \xinit \\
        &&& x^{k+1} = f(x^k, u^k) \quad\quad \raisebox{0.3ex}{$\scriptstyle(0 \le k \lt N)$}
    \end{aligned}
\end{equation}
The function \(f : \R^\nnx \times \R^\nnu \to \R^\nnx\) models the discrete-time,
nonlinear dynamics of the system, which starts from an initial state \(\xinit\). 
The functions \(h_k : \R^\nnx \times \R^\nnu \to \R^{n_y^k}\) for \(0 \le k \lt N\)
and \(h_N : \R^\nnx \to \R^{n_y^N}\) can be used to represent the (possibly time-varying) output mapping of the
system,
and the convex functions \(\ell_k : \R^{n_y^k} \to \R\) and
\(\ell_N : \R^{n_y^N} \to \R\) define the stage costs and the terminal cost
respectively.

The problem \eqref{eq:OCP} can be transformed into formulation \eqref{eq:P} as follows.
Recursively define the state transition function \(\Phi^k\) as \(\Phi^0(u) \defeq \xinit\) 
and \(\Phi^{k+1}(u) \defeq f\big(\Phi^k(u), u^k\big)\). Define \(G\) as the function
that maps a sequence of inputs to the interleaved states and inputs over the horizon, 
\(G(u) = \begin{pmatrix}
    \Phi^0(u), & u_0, & \Phi^1(u), & u_1, & \dots, & \Phi^N(u)
\end{pmatrix}\).
Using this definition, the \textit{single-shooting} or \textit{sequential} formulation of problem \eqref{eq:OCP}
is an instance of \eqref{eq:P}, with \(\ell = \ell_0 \oplus \dots \oplus \ell_N\),
\(h = h_0 \times \dots \times h_N\), \(F = h \circ G\), \(\psi = \ell \circ F\) and \(g = \delta_{\U}\).
Specifically,
\begin{equation}\label{eq:SS-OCP} \tag{SS-OCP}
    \begin{aligned}
        &\minimize_u && \ell\big(h\big(G(u)\big)\big) \\
        &\subjto && u \in \U. \\
    \end{aligned}
\end{equation}

\subsection{Gauss--Newton approximations for optimal control} \label{sec:gn-for-oc}
By specializing the Gauss--Newton QP \eqref{eq:eqp} for this class of optimal
control problems,
and by exploiting the separable structure of the objective function,
the Gauss--Newton step can be shown to be the solution to
the equality-constrained, finite-horizon, linear quadratic regulator problem \eqref{eq:ocp-qp}.

For the sake of readability, we defined the following variables.
\vspace{0.5em}
\begin{equation}
    \makebox[0pt][c]{\(
    \begin{aligned}
    &\begin{aligned}
        \bar x^k &\defeq \Phi^k(\bar u) & \hhbar^k &\defeq h_k(\barxuk) \\
        A_k &\defeq \jac_f^x(\barxuk) &
        B_k &\defeq \jac_f^u(\barxuk) \\
        q^k &\defeq \tp{\jac_{h_k}^x\!(\barxuk)} \nabla \ell_k(\hhbar^k)\;\; &
        r^k &\defeq \tp{\jac_{h_k}^u\!(\barxuk)} \nabla \ell_k(\hhbar^k) \\
        \Lambda_k &\defeq \partial^2 \ell_k(\hhbar^k) \\
    \end{aligned} \\
    &\begin{aligned}
        Q_k &\defeq \tp{\jac_{h_k}^x\!(\barxuk)} \Lambda_k\, \jac_{h_k}^x\!(\barxuk) \\
        S_k &\defeq \tp{\jac_{h_k}^u\!(\barxuk)} \Lambda_k\, \jac_{h_k}^x\!(\barxuk) \\
        R_k &\defeq \tp{\jac_{h_k}^u\!(\barxuk)} \Lambda_k\, \jac_{h_k}^u\!(\barxuk) \\
    \end{aligned}
    \end{aligned}
    \)}
\end{equation}

In order to transform \eqref{eq:ocp-qp} into a standard linear quadratic
regulator formulation, eliminate the fixed variables \(u_{\K}\).
The result is the problem \eqref{eq:ocp-qp-elim}, where we used the following
definitions.
\begin{equation}
    \begin{aligned}
        \hat S_k &\defeq S_{k}{\scriptstyle[\J\!,\,\cdot\,]} &
        \hat R_k &\defeq R_{k}{\scriptstyle[\J\!,\J]} \\
        \hat q_k &\defeq q^k + \tp S_{k}{\scriptstyle[\,\cdot\,,\K]}\, u^k_{\K} &\hspace{1.6em}
        \hat r_k &\defeq r^k_{\!\J} + R_{k}{\scriptstyle[\J\!,\K]}\, u^k_{\K} \\
        \hat B_k &\defeq B_{k}{\scriptstyle[\,\cdot\,,\J]} &
        \hat c_k &\defeq B_{k}{\scriptstyle[\,\cdot\,,\K]}\, u^k_{\K} \\
    \end{aligned}
\end{equation}

\begin{rem}
    In the absence of box constraints, we have \(\K = \emptyset\), and the algorithm
    reduces to the iterative linear quadratic regulator (ILQR) method for nonlinear 
    MPC of \cite[]{torodov_ilqr} with a line search.
\end{rem}

\subsection{Handling state constraints}\label{sec:state-constr}

Consider a standard state-constrained finite-horizon optimal control problem of
the following form.
\begin{equation}\label{eq:OCP-state-constr} \tag{SC-OCP}\hspace{-1em}
    \begin{aligned}
         & \minimize_{u,x} &  & \tfrac12\sum_{k=0}^{N-1} \bigg[ \normsq{x^k - \xref}_Q + \normsq{u^k -\uref}_R \bigg] \hspace{-1.1em} \\
         &                 &  & \hspace{-1.21em}+ \tfrac12 \normsq{x^N - \xref}_{Q_N}                                                 \\
         & \subjto         &  & u \in \U                                                                                              \\
         &                 &  & \begin{aligned}
             & x^0 = \xinit                                                   \\
             & x^{k+1} = f(x^k, u^k) &  & \quad \raisebox{0.3ex}{$\scriptstyle(0 \le k \lt N)$} \\
             & c_k(x^k) \in \D_k       &  & \quad \raisebox{0.3ex}{$\scriptstyle(0 \le k \le N)$} \\
        \end{aligned}
    \end{aligned}
\end{equation}
As before, \(f\) describes the possibly nonlinear discrete-time dynamics,
\(\xinit\) is the initial state of the system, \(\xref\) is the reference state,
and \(\uref\) the reference input. The inputs are constrained by the box \(\U\),
and some smooth, possibly nonlinear function \(c_k\) of the states enables the
representation of general equality and inequality constraints by constraining
its image to the box \(\D\).

It is common practice to relax the state constraints by means of a penalty
method. That is, the hard constraints are turned into soft constraints
by adding them as quadratic penalty terms to the objective function,
e.g. \(\frac{\mu}{2}\dist_{\D_k}^2\!\!\big(c_k(x^k)\big)\) for some sufficiently
large \(\mu > 0\).

Such a soft-constrained optimal control problem fits into the framework of
\eqref{eq:SS-OCP} by defining
\begin{equation}
    \begin{aligned}
        \ell_k(x, u, z) & \defeq \tfrac12 \normsq{x-\xref}_Q
        + \tfrac12 \normsq{u-\uref}_R + \tfrac{\mu_k}2 \dist^2_{\D_k}(z),\hspace{-1.5em} \\
        \ell_N(x, z)    & \defeq \tfrac12 \normsq{x-\xref}_{Q_N}
        + \tfrac{\mu_N}2 \dist^2_{\D_N}(z),                               \\
        h_k(x, u)       & \defeq \begin{pmatrix}
            x, &
            u, &
            c_k(x) \vphantom{\big|}
        \end{pmatrix},        \\
        h_N(x)          & \defeq \begin{pmatrix}
            x, &
            c_N(x) \vphantom{\big|}
        \end{pmatrix}.        \\
    \end{aligned}
\end{equation}

Because of the squared distance, the cost \(\ell\) is no longer twice
differentiable, but its gradient \(\nabla\ell\) is locally Lipschitz continuous,
and hence its Clarke generalized Jacobian \(\partial^2\ell\) is well defined and nonempty \cite[Prop.~7.1.4]{pang}.
Additionally, the gradient is semismooth, so Proposition~\ref{prop:lna} applies.

The following proposition gives a sufficient condition for the solution to 
the Gauss--Newton QP \eqref{eq:eqp} to be uniquely defined.
\begin{prop}
    If the cost matrix \(R\) is positive definite, \(Q\) is positive semidefinite,
    and \(\mu_k \ge 0\) for all \(k\),
    then the Gauss--Newton matrix \(\Hgn\) for the soft-constrained optimal
    control problem is positive definite.
\end{prop}
\begin{pf}By algebraic manipulations of \(\Hgn\).\\
    Because of the block-diagonal structure of \(\partial^2\ell\) and \(\jac_h\),
    their product \(L \defeq \tp{\jac_h}\partial^2\ell\,\jac_h\) is also block-diagonal,
    with blocks of the form 
    \begin{equation*}
        \begin{pmatrix}
            Q + \tp C_k M_k C_k & 0 \\ 0 & R
        \end{pmatrix} \possdef 0,
    \end{equation*}
    where \(C_k \defeq \jac_{c_k}(x^k)\) and \(M_k \in \partial^2\big(\frac{\mu}{2}\dist_{\D_k}^2\!(c_k(x^k))\big)\). Because of the
    structure of \(G\) (it includes the identity map of \(u\)), the block
    rows of \(\jac_G(u)\) that correspond to the inputs have full rank (they contain \(\nnu\times\nnu\) identity matrices)
    and line up with the positive definite blocks \(R\) in \(L\).
    Hence, the full product
    \(\Hgn = \tp{\jac_G(u)}L\,\jac_G(u)\) is positive definite.
    \hfill\(\square \)
\end{pf}

\section{Algorithmic details}\label{sec:algorithm}
We will now explore algorithms for efficiently solving
\eqref{eq:ocp-qp-elim} to obtain the Gauss--Newton step \(\delu\) that can
be used to accelerate PANOC.

For the sake of self-containedness, the PANOC$^+$ method from \cite[]{de_marchi_proximal_2022}
is given in Algorithm~\ref{alg:panoc}. It has been specialized to use
the Gauss--Newton step \(\delu\) derived in Section~\ref{sec:gn-panoc}. Unlike
the original version of PANOC$^+$ with an L--BFGS accelerator, a Gauss--Newton
step can be computed from the very first iteration.

\begin{algorithm2e}
    \newcommand\assign{\leftarrow}
    \newcommand\nupo{\nu\hspace{-0.6pt}+\hspace{-0.8pt}1}
    \newcommand\numo{\nu\hspace{-0.4pt}-\hspace{-0.8pt}1}
    \DontPrintSemicolon
    \KwIn{initial guess $u^{(0)}$,
    initial step size $\gamma_0 > 0$,
    parameters $\alpha, \beta \in (0, 1)$}
    \KwOut{$u^\star$}
    \vspace{0.5em}
    $\hat u^{(0)} \assign T_{\gamma_{0}}(u^{(0)}), \quad p^{(0)} \assign \hat u^{(0)} - u^{(0)} $ \;
    $\nu \assign 1$\;
    \While{\normalfont Stopping criterion not satisfied for $u^{(\numo)}$}
    { \vspace{0.2em}
        Compute $\delu$ from \eqref{eq:eqp} with \(\bar u \defeq u^{(\numo)}\)\;
        $\gamma_{\nu} \assign \gamma_{\numo}, \quad \tau \assign 1$ \;
        \vspace{0.6ex}
        \hspace{-0.6em}\raisebox{2.5ex}{\makebox[0pt][l]{\(\triangleright\)}}\hspace{0.6em}%
        $u^{(\nu)} \assign u^{(\numo)} + (1-\tau)\, p^{(\numo)} + \tau\, \delu$\;
        $\hat u^{(\nu)} \assign T_{\gamma_{\nu}}(u^{(\nu)}), \quad p^{(\nu)} \assign \hat u^{(\nu)} - u^{(\nu)} $ \;
        \If{$\psi(\hat u^{(\nu)}) > \psi(u^{(\nu)}) + \ttp{\nabla \psi(u^{(\nu)})} p^{(\nu)} + \tfrac{\alpha}{2\gamma_{\nu}} \normsq{p^{(\nu)}} $}
        {
            $\gamma_{\nu} \assign \gamma_{\nu} / 2, \quad \tau \assign 1$ and go to \(\triangleright\)\;
        }
        \If{$\fbe[\gamma_{\nu}](u^{(\nu)}) > \fbe[\gamma_{\numo}](u^{(\numo)}) - \beta \tfrac{1-\alpha}{2\gamma_{\numo}} \normsq{p^{(\numo)}}$}{
            $\tau \assign \tau / 2$ and go to \(\triangleright\)\;
        }
        $\nu \assign \nu + 1$\;
    }
    $u^\star \assign T_{\gamma_{\numo}}(u^{(\numo)})$\;
    \caption{PANOC$^+$ \cite[Algorithm 2]{de_marchi_proximal_2022} with Gauss--Newton acceleration}\label{alg:panoc}
\end{algorithm2e}

\subsection{Evaluation of the objective and its gradient}
Application of PANOC to problem \eqref{eq:SS-OCP} requires efficient evaluation
of the cost function \(\psi = \ell \circ h \circ G\) and its gradient. This can
be achieved by performing a
forward simulation (Algorithm~\ref{alg:lqr-sim}) followed by a backward sweep
(Algorithm~\ref{alg:lqr-grad}). The backward sweep only requires the evaluation
of gradient-vector products, but the Jacobian matrices \(A_k\) and \(B_k\) of
the dynamics can later be reused for the computation of the Gauss--Newton step.
\begin{algorithm2e}
    \newcommand\assign{\leftarrow}
    \DontPrintSemicolon
    \KwIn{$\bar u, \xinit$}
    \KwOut{$\psi, \bar x, \hhbar$}
    $\bar x^0 \assign \xinit$\;
    $\psi \assign 0$\;
    \For{$k = 0,\, ...,\, N-1$}
    {
        $\bar x^{k+1} \assign f(\barxuk)$\;
        $\hhbar^k \assign h_k(\barxuk)$\;
        $\psi \assign \psi + \ell_k(\hhbar^k)$\;
    }
    $\hhbar^N \assign h_N(\bar x^N)$\;
    $\psi \assign \psi + \ell_N(\hhbar^N)$\;
    \caption{Forward simulation}\label{alg:lqr-sim}
\end{algorithm2e}
\begin{algorithm2e}
    \newcommand\assign{\leftarrow}
    \DontPrintSemicolon
    \KwIn{$\bar u^k, \bar x^k, \hhbar^k$}
    \KwOut{$\nabla\psi, A_k, B_k, q^k, r^k$}
    $\lambda^N \assign \tp{\jac_{h_N}\!(\bar x^N)} \nabla \ell_N(\hhbar^N) $\;
    \For{$k = N - 1,\, ...,\, 0$}
    {
        $\begin{pmatrix}
            A_k & B_k
        \end{pmatrix} \assign \jacf(\barxuk)$\;
        $q^k \assign \tp{\jac_{h_k}^x\!(\barxuk)} \nabla \ell_{k}(\hhbar^k)$\;
        $r^k \assign \tp{\jac_{h_k}^u\!(\barxuk)} \nabla \ell_{k}(\hhbar^k)$\;
        $\nabla_{\!u^{\hspace{-0.6pt}k}}\psi \assign r^k + \tp B_k \lambda^{k+1}$\;
        $\lambda^k \assign q^k + \tp A_k \lambda^{k+1}$\;
    }
    \caption{Backward gradient evaluation}\label{alg:lqr-grad}
\end{algorithm2e}

\subsection{Solution of the LQR problem}
The Gauss--Newton step \(\delu\) can be computed as the solution to
\eqref{eq:ocp-qp-elim} using LQR factorization and LQR solution
routines based on the Riccati recursion \cite[\S 8.8.3]{rawlings2017model},
\cite[Alg.~3-4]{patrinoslqrsolve}.
These routines, specialized to the problem at hand, are listed in
Algorithms~\ref{alg:lqr-factor} and \ref{alg:lqr-solve}.

\begin{algorithm2e}[h]
    \newcommand\assign{\leftarrow}
    \DontPrintSemicolon
    \KwIn{$Q_k, \hat S_k, \hat R_k, \hat q_k, \hat r_k, A_k, \hat B_k, \hat c_k$}
    \KwOut{$K_k, e_k$}
    $P_N \assign Q_N$\;
    $s_N \assign \hat q_N$\;
    \For{$k = N - 1,\, ...,\, 0$}
    {
        $\bar R \assign \hat R_k + \tp{\hat B_k} P_{k+1} \hat B_k$\;
        $\bar S \assign \hat S_k + \tp{\hat B_k} P_{k+1} A_k$\;
        $y \assign P_{k+1} \hat c_k + s_{k+1}$\;
        $K_k \assign -\inv{\bar R} \bar S$\;
        $e_k \assign -\inv{\bar R} (\tp{\hat B_k} y + \hat r_k)$\;
        $s_k \assign \tp{\bar S} e_k + \tp A_k y + \hat q_k$\;
        $P_k \assign Q_k + \tp A_k P_{k+1} A_k + \tp{\bar S}\! K_k$\;
    }
    \caption{LQR factor}\label{alg:lqr-factor}
\end{algorithm2e}
\begin{algorithm2e}[h]
    \newcommand\assign{\leftarrow}
    \DontPrintSemicolon
    \KwIn{$A_k, B_k, K_k, e_k, \delu_\K$}
    \KwOut{$\delu_{\!\J}, \delx$}
    $\delx^0 \assign 0$\;
    \For{$k = 0,\, ...,\, N - 1$}
    {
        $\delu^k_{\!\J} \assign K_k \delx^k + e_k$\;
        $\delx^{k+1} \assign A_k \delx^k + B_k \delu^k$\;
    }
    \caption{LQR solve}\label{alg:lqr-solve}
\end{algorithm2e}

An important observation is that the cost for the computation of the
Gauss--Newton direction using these routines scales \textit{linearly} with the
horizon length \(N\). In the worst case, when \(\K(\bar u) = \emptyset\),
Algorithm~\ref{alg:lqr-factor} requires the factorization of \(N\) matrices of
size \(\nnu \times \nnu\) and some matrix products.
In contrast, general direct solution methods for
system \eqref{eq:struc-panoc-sys} require a single factorization of a much
larger \(\nnu N \times \nnu N\) matrix, with a cost that scales
\textit{cubically} with \(N\).

\subsection{Practical considerations} \label{sec:practical-considerations}
For iterates that are far from the solution, the quadratic Gauss--Newton model
might not approximate the actual function well, and the Gauss--Newton step might
not perform much better than an L--BFGS step.
Considering the significant difference in computational cost between
Gauss--Newton and L--BFGS (the former requires evaluation of the Jacobians of
the dynamics, matrix factorizations and multiplications, whereas the latter
only requires a limited number of vector operations),
we propose to only compute the Gauss--Newton step every \(k_\mathrm{GN} \ge 1\)
iterations. In between, much cheaper structured \panoc{} L--BFGS steps
are used \cite[\S III]{pas2022alpaqa}. When eventually a Gauss--Newton step is
accepted by the line search with step size \(\tau = 1\),
the algorithm continues to perform Gauss--Newton steps, for as long as they
keep getting accepted with unit step size. 
Using this technique, the algorithm initially maintains a relatively low
cost per iteration, and eventually enjoys the fast local convergence of the
more expensive Gauss--Newton steps. This will be corroborated experimentally
in the following section.

\section{Experimental results}\label{sec:experiments}
In this section, the PANOC algorithm with Gauss--Newton acceleration is applied
to a nonlinear, input-constrained model predictive control problem, and its performance is
compared to the approximate structured PANOC algorithm with L--BFGS
acceleration from \cite[]{pas2022alpaqa}.
As a benchmark, we consider the optimal control of a \textit{``chain of masses connected by springs''}
described by \cite[]{chain}. 
One side of the chain is fixed, and the
other side is attached to an actuator. A disturbance is applied to the system,
and the goal of the controller is to bring the chain back to a steady state,
with the actuator at a predetermined target position. 
The input constraints limit the velocity of the actuator to 
\(1\,\mathrm{m}/\mathrm{s}\) along each axis. Unless specified otherwise, we use
the parameter values listed in \cite[]{chain}.

The software package CasADi \cite[]{Andersson2019} is used to model and
discretize the problem using a fourth-order Runge--Kutta integrator,
and the resulting subroutines
for evaluating the dynamics, the stage cost and terminal cost functions,
as well as their derivatives are compiled, and used in an optimized \Cpp{} implementation of
Algorithms~\ref{alg:panoc}--\ref{alg:lqr-solve}, based on \textsc{alpaqa} \cite[]{pas_alpaqa_github}.%
\footnote{The Python source code to reproduce the results in this section can be found at \href{https://github.com/kul-optec/panoc-gauss-newton-ifac-experiments}{{github.com/kul-optec/panoc-gauss-newton-ifac-experiments}}.
All experiments were carried out using an Intel Core i7-7700HQ CPU at 2.8 GHz.}

\subsection{Number of iterations}
\begin{figure}
    \centering
    \includegraphics[width=0.8\linewidth]{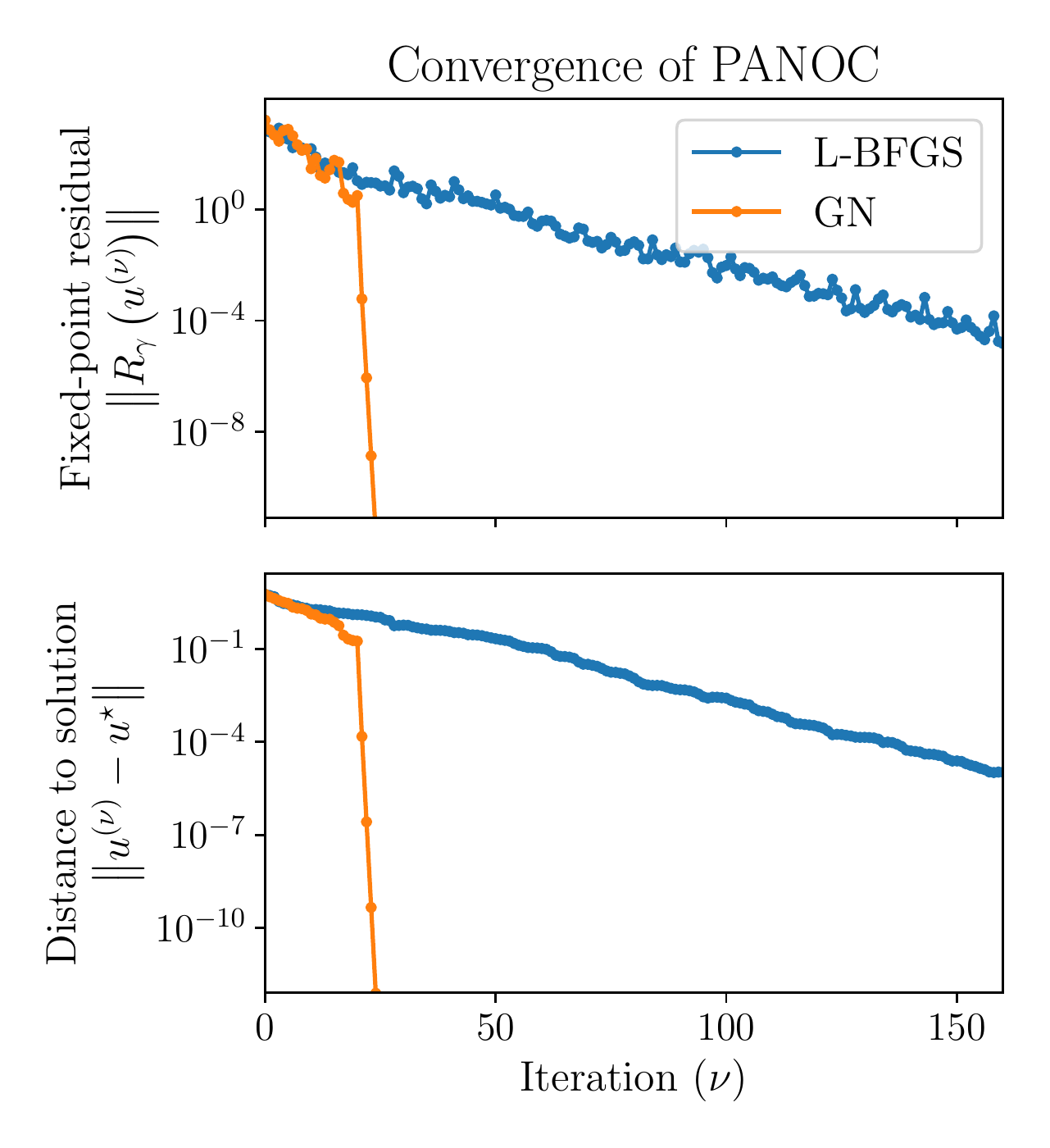}
    \vspace{-1em}
    \caption{Comparison of the convergence of structured \panoc{} with L--BFGS and \panoc{} with the proposed
    Gauss--Newton accelerator (\(k_\mathrm{GN} = 1\)), when applied to the chain of masses MPC benchmark.} \label{fig:conv-iter}
\end{figure}

In a first experiment, the convergence in terms of the number of iterations is compared for
the PANOC algorithm with Gauss--Newton acceleration as described in this
publication, and for the structured PANOC algorithm with L--BFGS acceleration without the
off-diagonal Hessian--vector term from \cite[]{pas2022alpaqa}. For the Gauss--Newton
accelerator, the parameter \(k_\mathrm{GN}\)
from Section~\ref{sec:practical-considerations} is set to one (i.e. a
Gauss--Newton step is computed on each PANOC iteration). The L--BFGS memory is set to 40,
equal to the length of the horizon.
Figure~\ref{fig:conv-iter} shows the convergence of the two algorithms. Initially,
they both perform similarly, but after around 20 iterations, the Gauss--Newton
directions are accepted with unit step size, enabling very fast linear
convergence.

It should be noted that similar graphs in terms of absolute solver run time
would look quite different: even though the reduction of the residual per 
iteration is comparable for the first 20 iterations,
the computational cost per iteration for the Gauss--Newton accelerator is around
one order of magnitude higher than for the L--BFGS accelerator. This can be
greatly improved by increasing \(k_\mathrm{GN}\).

\subsection{Run time in function of horizon length}
\begin{figure}
    \centering
    \includegraphics[width=\linewidth]{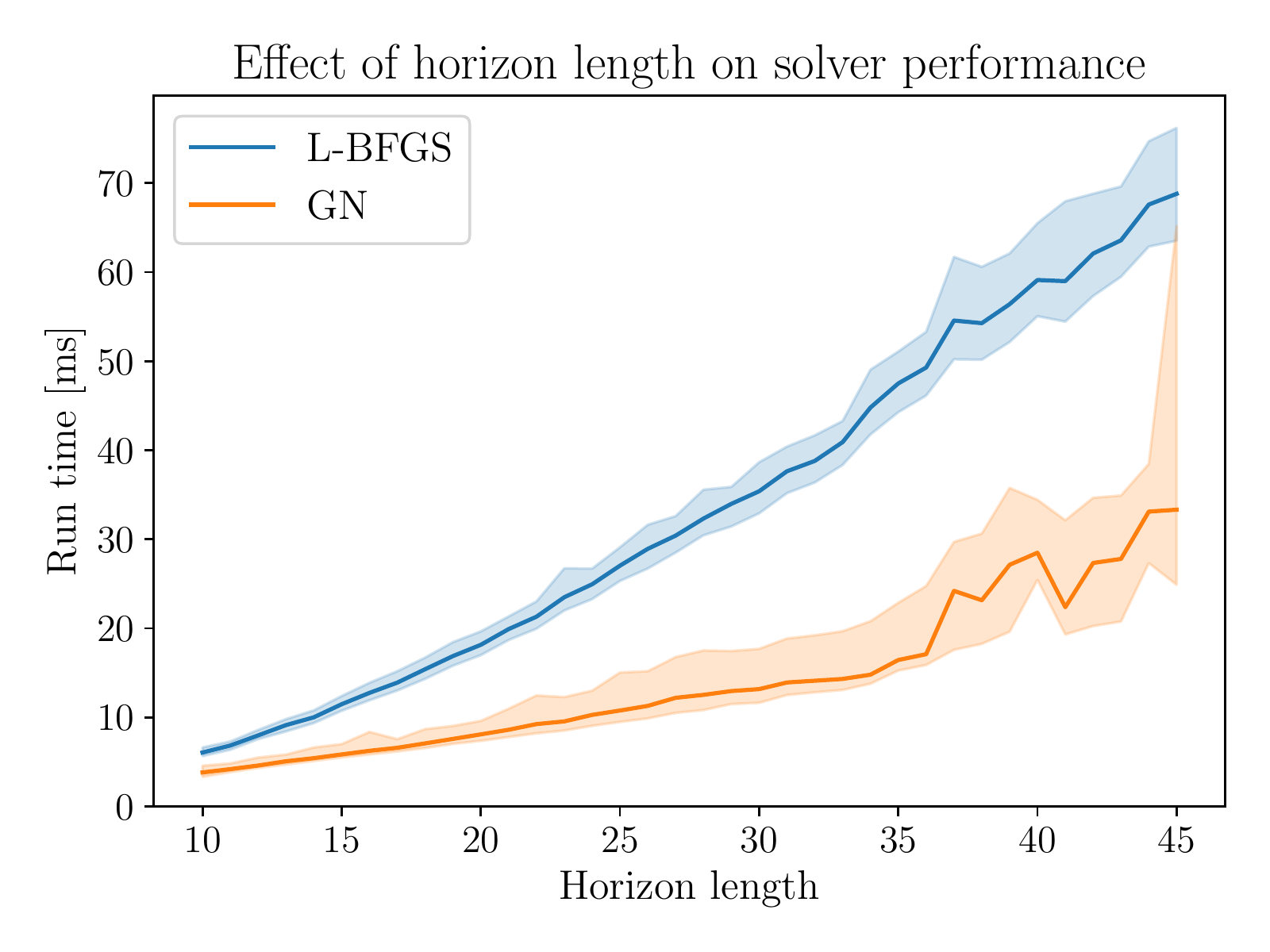}
    \vspace{-1em}
    \caption{Median solver run time over the 256 test problems for each horizon length, for structured \panoc{} with L--BFGS and \panoc{} with the Gauss--Newton accelerator (\(k_\mathrm{GN} = 30\)).
    The shaded area indicates the P10 and P90 percentiles.}\label{fig:horiz}
\end{figure}

In a second experiment, we explore the effect of the horizon length on the
solver run time. For each horizon length between \(N=10\) 
and \(N=45\), 256 optimal control problems are composed, each with a different
initial state \(\xinit\), generated by applying uniformly random inputs in 
\([-1, 1]\) for five time steps.
The parameter \(k_\mathrm{GN}\) described in Section~\ref{sec:practical-considerations}
was set to 30 for this experiment, and the L--BFGS memory was set equal to the
horizon length \(N\). The solvers declare convergence when 
\(\left\| u^{(\nu)} - \proj{\U}\!\left(u^{(\nu)} - \nabla \psi(u^{(\nu)})\right) \right\| \le 10^{-10}\).
The run times of both algorithms (structured PANOC with L--BFGS, and
PANOC with Gauss--Newton acceleration) are reported in Figure~\ref{fig:horiz}.
The algorithm with Gauss--Newton acceleration is more than twice as fast as
the L--BFGS variant, and the run time scales not much worse than linearly with
the horizon length \(N\), although longer horizons appear to be more challenging.

\subsection{Model predictive control}
\begin{figure}
    \centering
    \includegraphics[width=\linewidth]{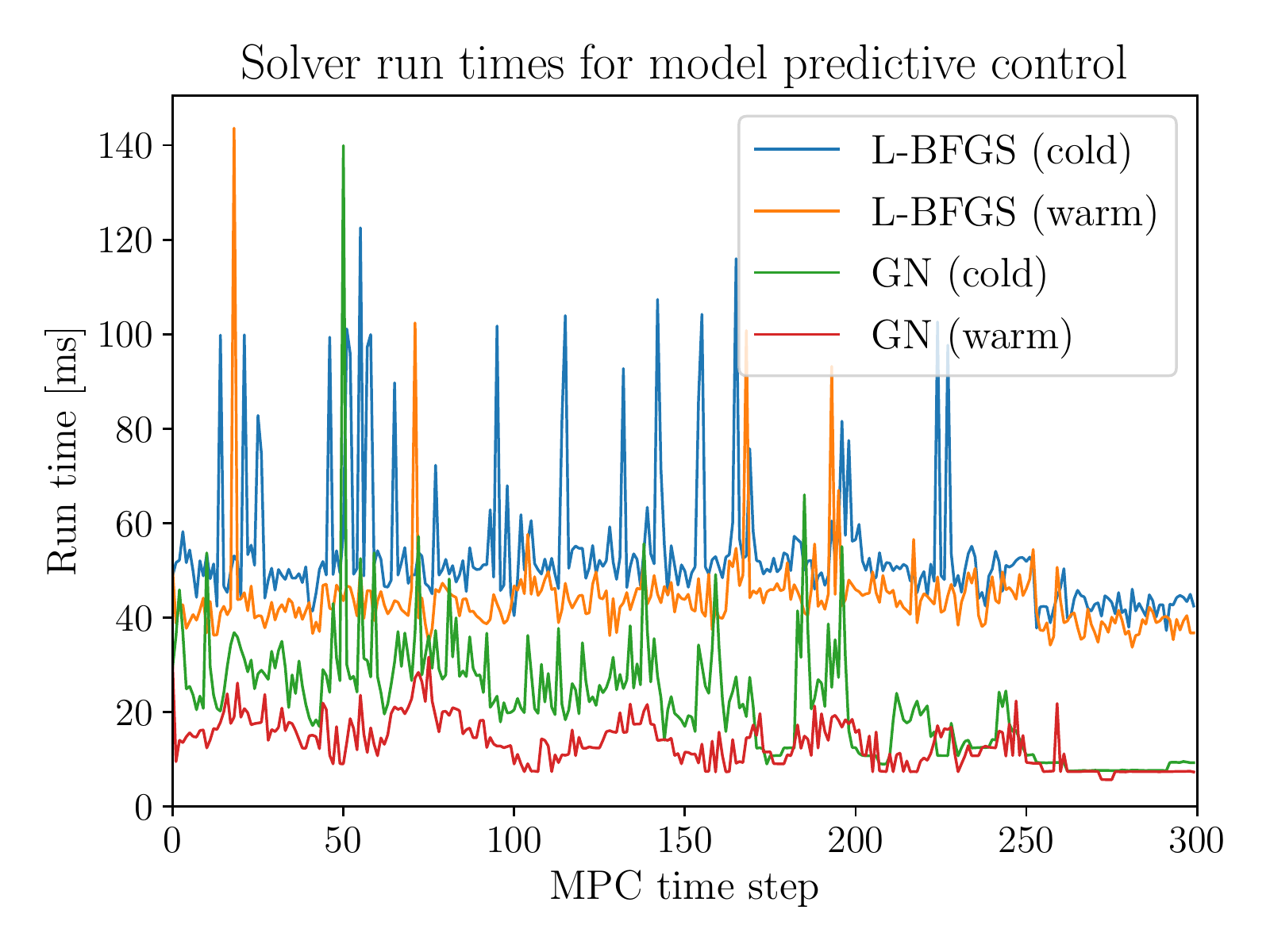}
    \vspace{-1em}
    \caption{Solver run times for structured \panoc{} with L--BFGS and \panoc{} with the Gauss--Newton accelerator (\(k_\mathrm{GN} = 10\)) when applied to a model predictive control problem.
    For data labeled \textit{warm}, the shifted solution of the previous time step is used as initial guess for the solvers, whereas it is set to zero for data labeled \textit{cold}.}\label{fig:mpc}
\end{figure}

Finally, both solvers are applied in a closed-loop controller. A disturbance 
of \([-1, 1, 1]\, \mathrm{m/s}\) is applied for five time steps, and the
system with the MPC controller is subsequently simulated for one minute. The run times of the
two solvers described earlier are reported in Figure~\ref{fig:mpc}.
The Gauss--Newton solver (with \(k_\mathrm{GN} = 10\)) outperforms the 
L--BFGS-based solver in terms of both average and worst-case run time.
The fast local convergence of Gauss--Newton is especially noticeable when the
initial guess is close to the solution, e.g. by warm starting the solver using
the shifted solution from the previous time step, and when the system starts to
settle near the end of the simulation. For reference,
the popular \texttt{Ipopt} solver \cite[]{ipopt} requires around 1.7 seconds to
solve the first OCP (invoked from CasADi, without just-in-time compilation),
which is over 50 times longer than the 30 ms required by the \panoc{} solver
with Gauss--Newton acceleration.

\section{Conclusion}\label{sec:conclusion}
In this paper, we extended the \panoc{} algorithm to enable acceleration using
Gauss--Newton directions. We showed how the structure of optimal control
problems can be exploited to efficiently compute these Gauss--Newton directions
using the Riccati recursion, in such a way that the computational cost scales
linearly with the horizon length.
Performance of the proposed methods was then compared to
a previous variant of \panoc: we reported a speedup by a factor of two for a
challenging optimal control benchmark problem.

An open-source \Cpp{} implementation of the algorithm is under active development
in the \textsc{alpaqa} GitHub repository \cite[]{pas_alpaqa_github}.
Using the techniques outlined in Section~\ref{sec:state-constr}, the method can
be integrated into \textsc{alpaqa}'s augmented Lagrangian and quadratic penalty
framework.
Further performance improvements could be achieved by exploiting the
sparsity of the Jacobians \(A_k\) and \(B_k\) and/or by employing specially tailored
linear algebra routines such as BLASFEO \cite[]{Frison2018BLASFEOBL}.

\bibliography{paper}

\end{document}